\documentclass{amsproc}
\usepackage{amsfonts,amssymb,amsmath,amsthm,amscd,mathtools,multicol,tikz, tikz-cd,caption,enumerate,mathrsfs,thmtools,cite}
\usepackage[pagebackref=true, colorlinks, linkcolor=blue, citecolor=red]{hyperref}
\usepackage{fullpage}

\usepackage{palatino}
\makeatletter
\makeindex

\declaretheoremstyle[headfont=\normalfont]{normalhead}
\newtheorem{theorem}{Theorem}[section]
\newtheorem{proposition}[theorem]{Proposition}
\newtheorem{corollary}[theorem]{Corollary}
\newtheorem{definition}[theorem]{Definition}
\newtheorem{lemma}[theorem]{Lemma}
\theoremstyle{remark}
\newtheorem{remark}[theorem]{\bf{Remark}}
\newtheorem{example}[theorem]{\bf{Example}}

\numberwithin{equation}{section}

\newcommand{\bt}{\begin{theorem}}
\newcommand{\et}{\end{theorem}}
\newcommand{\bco}{\begin{corollary}}
\newcommand{\eco}{\end{corollary}}
\newcommand{\bd}{\begin{definition}}
\newcommand{\ed}{\end{definition}}
\newcommand{\bp}{\begin{problem}}
\newcommand{\ep}{\end{problem}}
\newcommand{\bl}{\begin{lemma}}
\newcommand{\el}{\end{lemma}}
\newcommand{\bprop}{\begin{proposition}}
\newcommand{\eprop}{\end{proposition}}
\newcommand{\br}{\begin{remark}}
\newcommand{\er}{\end{remark}}
\newcommand{\bpf}{\begin{proof}}
\newcommand{\epf}{\end{proof}}
\newcommand{\bex}{\begin{example}}
\newcommand{\eex}{\end{example}}

\newcommand{\C}{\mathbb{C}}

\newcommand{\pv}{E}

\newcommand{\R}{\overline{R}}
\newcommand{\Sm}{\overline{S}}
\newcommand{\fieldfrac}{\mathcal Q}
\newcommand{\variety}{\mathscr{V}}
\newcommand{\gal}{\mathscr{G}}
\newcommand{\galsub}{\mathscr{H}}
\newcommand{\uni}{\mathscr{U}}
\newcommand{\borel}{\mathscr{B}}

\newcommand{\id}{\mathfrak{a}}
\newcommand{\ib}{\mathfrak{b}}
\newcommand{\p}{\mathfrak{p}}
\newcommand{\q}{\mathfrak{q}}

\renewcommand{\a}{\alpha}
\newcommand{\s}{\sigma}
\newcommand{\Lde}{\mathscr L}

\begin{document}

\title{A NOTE ON LIOUVILLIAN PICARD-VESSIOT EXTENSIONS}
\author{Ursashi Roy \ \ and \ \ Varadharaj R.  Srinivasan}
\address{Indian Institute of Science Education and Research, Sector 81, Mohali 140306, India}


\maketitle

\begin{abstract}In this paper, we prove a new characterization theorem for  Picard-Vessiot extensions whose  differential Galois groups have solvable identity components.\end{abstract}

\section{INTRODUCTION}\label{intro} Throughout this article, we fix a differential field\footnote{A field equipped with a single derivation map denoted by $'$.} $F$ of characteristic zero with an algebraically closed field of constants $C:=\{x\in F\  | \ x'=0\}.$  Let $E$ be a Picard-Vessiot extension of $F,$ $K$ be a differential field intermediate to $E$ and $F$ and $T(K|F)$ be the set of all solutions in $K$ of all linear homogeneous differential equations over $F.$ It is known that  $T(E|F)$  is a finitely generated simple differential  $F-$algebra whose field of fractions $\fieldfrac(T(E|F))$ equals the differential field $E$. However, if $r\in E$ and  $f,g\in T(E|F)$ are elements such that $r=f/g$ then it is possible that neither $f$ nor $g$ belong to the differential field $\langle r\rangle,$ generated by $F$ and $r.$ Thus $\fieldfrac(T(K|F))$ could be a proper subfield of $K$. For example, consider the ordinary differential field $(\C(x),')$ of complex rational functions with derivation $':=d/dx.$ Let $E$ be a Picard-Vessiot extension of the Airy differential equation $\Lde(Y):=Y''-xY=0.$ Then the differential Galois group is isomorphic to $SL(2,\C)$ as algebraic groups. The differential field $K$ fixed by the subgroup of upper triangular matrices in $SL(2,\C)$ is of the form  $K=\C(x)(w),$ where $w$ is transcendental over $\C(x)$ and $w$ is a solution of the Ricatti equation $w'=x-w^2$ and that $T(K|\C(x))=T(E|\C(x))\cap K=\C(x)$ (\cite{Magid1999}, pp. 86-87). Therefore, it is natural to ask for a characterization theorem of those Picard-Vessiot extensions whose intermediate differential fields $K$ are the field of fractions of $T(K|F).$  
 
A differential field extension $E$ of $F$ is called a \emph{liouvillian Picard-Vessiot} extension if $E$ a liouvillian extension as well as a Picard-Vessiot extension of $F.$ Liouvillian Picard-Vessiot extensions are characterized by their differential Galois groups having a solvable identity component (\cite{Put-Singer2003}, Theorem 1.43). Using this fact and the well-known structure theorem of $T(E|F)$ (\cite{Magid1999}, Theorem 5.12), we prove  that if $E$ is a Picard-Vessiot extension of $F$ then $\fieldfrac(T(K|F))=K$ for every differential field intermediate to $E$ and $F$ if and only if $E$ is a liouvillian Picard-Vessiot extension of $F;$ in which case, we also show that $T(K|F)$ is a finitely generated simple differential $F-$algebra (Theorem \ref{NSliouvillian}). If the differential Galois group of a liouvillian Picard-Vessiot extension $E$ of $F$ is connected then given any intermediate differential subfield $K$, we find a tower of differential fields $$F=K_0\subseteq K_1\subseteq\cdots\subseteq K_{m-1}\subseteq K_m=K,$$ such that  for each $1\leq i\leq m$, there is a $t_i\in T(E|F)$ such that  $K_i=K_{i-1}(t_i)$  and  $t'_i=a_it_i+b_i$ for some $a_i\in F$ and $b_i\in K_{i-1}$ (Corollary \ref{structureCLPV}). 

We also prove that if $E$ is an arbitrary Picard-Vessiot extension of $F$ and $K$ is an intermediate differential field then $K$ contains a  finitely generated simple differential $F-$subalgebra $R$ such that  $\mathcal Q(R)=K$ (Theorem \ref{IDSS}). A structure theorem, similar to Corollary \ref{structureCLPV},  for relatively algebraically closed intermediate differential subfields of liouvillian extensions  can be found in \cite{Ravis2020}.
For fundamental results on Picard-Vessiot theory, we refer the reader to \cite{Magid1999} and \cite{Put-Singer2003}.

\section{PRELIMINARIES}\label{prelims}

In this section we record few definitions and results from Picard-Vessiot theory that are used in our proofs.  Let $F[\partial]$ be the ring of differential operators over $F$ and  $\Lde\in F[\partial]$ be a monic operator of order $n$. A \emph{Picard-Vessiot extension} $E$ of $F$ for $\Lde$ is a differential field extension of $F$ having the same field of constants as $F$ and satisfying the following conditions: \begin{enumerate}[(a)] \item The $C-$vector space $V$ of all solutions of $\Lde(Y)=0$ in $E$ is of dimension $n.$\\ \item  $E=F\langle V\rangle,$ that is, the smallest differential field containing $F$ and $V$ is $E$.
\end{enumerate} The differential Galois group, denoted by $\gal(E|F)$, is the group of all field automorphisms of $E$  that fixes the elements of $F$ and commutes with the derivation of $E.$   The differential Galois group stabilizes $V$ and thus it acts as a group of $C-$vector space (differential) automorphisms of $V.$ Since $E=F\langle V\rangle,$ the induced map $\phi: \gal(E|F) \to GL(V)$ is a faithful representation of groups. In fact, $\phi(\gal(E|F))$ can be shown to be a Zariski closed subgroup of $GL(V)$ and in this sense, $\gal(E|F)$ is seen as a linear algebraic group. The fundamental theorem of Picard-Vessiot theory provides a bijective correspondence between differential subfields intermediate to $E$ and $F$ and the Zariski closed subgroups of $\gal(E|F)$ in a way that closely resembles the polynomial Galois theory. If $\galsub$ is a closed subgroup of $\gal(E|F)$ and $K$ is an intermediate differential field then the bijective correspondence is given by the maps \begin{align*}K \longmapsto & \gal(E|K):=\{\s\in \gal(E|F)\ | \ \s(u)=u\ \text{for all}\ u\in K\}\\  \galsub\longmapsto & E^\galsub:=\{u\in E\ | \ \s(u)=u\ \text{for all}\ \s\in \galsub\}.\end{align*}
 The field fixed by $\gal(E|F)$ is $F;$ that is $E^{\gal(E|F)}=F.$ Let $K$ be a differential field intermediate to $E$ and $F.$ Then $K$ is a Picard-Vessiot extension of $F$ if and only if $\gal(E|K)$ is a closed normal subgroup of $\gal(E|F)$ and in which case, the differential Galois group $\gal(K|F)$ is isomorphic to the quotient group $\gal(E|F)/\gal(E|K).$  If an intermediate differential field $K$ is stabilized by the differential Galois group then $\gal(E|K)$ is a normal subgroup of $\gal(E|F)$ and consequently, $K$ is a Picard-Vessiot extension of $F.$  The algebraic closure of $F$ in $E$  is a finite Galois extension, which we denote by $F(x).$ Clearly, $F(x)$ is stabilized by $\gal(E|F)$ and in fact, $F(x)=E^{\gal(E|F)^0},$  where $\gal(E|F)^0$ is the connected component of $\gal(E|F).$ The quotient group $\gal(E|F)/\gal(E|F)^0$ coincides with the ordinary Galois group of $F(x)$ over $F.$

 The differential $F-$algebra $T(E|F),$ consisting of all solutions in $E$ of linear homogeneous differential equations over $F$, plays a very important role in Picard-Vessiot theory and it is well understood.  The following facts  on $T(E|F)$ are well-known: The differential Galois group stabilizes $T(E|F),$  $\fieldfrac(T(E|F))=E$  and if the $\gal(E|F)$ orbit set of an element $y\in E$ spans a finite dimensional $C-$vector space  then $y\in T(E|F).$ There is a structure theorem that describes $T(E|F)$ in terms of the  coordinate ring of $\gal(E|F)$ (\cite {Magid1999}, Theorem $5.12$): If $\overline{F}$ is an algebraic closure of $F$ then there is an $\overline{F}-$algebra isomorphism $$\overline{F}\otimes_F T(E/F)\longrightarrow \overline{F}\otimes_C C[\gal(E|F)].$$ Furthermore, the above isomorphism respects the $\gal(E|F)$ action. Here $\gal(E|F)$ acts trivially on $\overline{F}$ and acts as right translations on the coordinate ring $C[\gal(E|F)]$ of $\gal(E|F).$ When $\gal(E|F)$ is  a connected solvable group, it is also known that  \begin{equation} \label{connectedST} T(E|F)\simeq F\otimes_C C[\gal(E|F)],\end{equation} where again, the isomorphism is compatible with the action of $\gal(E|F)$(\cite{Magid1999}, Corollary $5.29$).

In this article we will be studying the $F-$algebra $T(K|F),$ where $K$ is an intermediate differential field of a Picard-Vessiot extension of $F.$ In view of the fundamental theorem, if $\galsub$ is a closed subgroup then $\galsub=\gal(E|K)$ for some intermediate differential field $K$ and  we have $T(K|F)=T(E|F)\cap K=T(E|F)^\galsub.$ As noted in the introduction section, it can happen that $T(K|F)=F.$  The characterization theorem we prove in this article says that $\gal(E|F)^0$ is solvable if and only if  $T(K|F)$ has "enough elements" in the sense that  $\fieldfrac(T(K|F))=K$ for every intermediate differential field $K.$ For the proof of our theorem, we will  rely on the structure theorem described in Equation \ref{connectedST}  along with the following proposition.

\bprop\label{intext} Let $F$ be a differential field of characteristic zero with an algebraically closed field of constants. Let $E$ be a Picard-Vessiot extension of $F$ and $F(x)$ be the algebraic closure of $F$ in $E.$ Let $K$ be a differential field intermediate to $F$ and $E.$ Then \begin{enumerate}[(a)]\item   $T(K(x)|F(x))=T(K(x)|F).$\\ \item $T(K(x)|F)$ is an integral extension of $T(K|F).$
\end{enumerate}\eprop

\bpf
 Every differential equation over $F$ is also a differential equation  $F(x)$ and thus it is clear that $T(K(x)|F)\subseteq T(K(x)|F(x)).$ Since $F(x)$ is finite dimensional $F-$vector space, for any $y\in F(x),$ there must be a nonnegative integer $m$ such that $y, y',\cdots, y^{(m)}$ are $F-$linearly dependent. Therefore, $F(x)\subseteq T(K(x)|F).$ Now let $y\in T(K(x)|F(x))\setminus F(x)$ and $\Lde=\partial^{(n)}+a_{n-1}\partial^{(n-1)}+\cdots+a_0\in F(x)[\partial]$ be a monic operator of order $n\geq 1$ such that $\Lde(y)=0.$   Let $V$ be the  set of all solutions of $\Lde$ in $E$ and for any $\s\in \gal(E|F),$ let $V_\s$ be the set of all solutions of $\Lde_\s=\partial^n+\s(a_{n-1})\partial^{(n-1)}+\cdots+\s(a_0)$ in $E.$ Observe that  $\s(V)=V_\s.$ Since $a_i\in F(x)$ for each $i$ and $E^{\gal(E|F)^0}=F(x),$ the orbit set of $a_i$ under the action of $\gal(E|F)$ is a finite set for each $i$. Therefore, there are only finitely many $\Lde_\s.$ Let $\s_0\in \gal(E|F)$ be the identity and $\Lde=\Lde_{\s_0}, \Lde_{\s_1},\cdots, \Lde_{\s_l}$ be the distinct operators. Let $W=V_{\s_0}+V_{\s_1}+\cdots+V_{\s_l}.$ Clearly, $W$ is a finite dimensional $C-$vectorspace. For any $\s\in \gal(E|F)$ and $y\in V_{\s_i},$ we have $\s(y)\in V_{\s\s_i}=V_{\s_j}\subseteq W.$ This implies that $W$ is also a $\gal(E|F)-$module. Thus, any $y\in W$ must be a solution of some operator in $F[\partial].$  That is, $y\in T(K(x)|F).$    Hence  \begin{equation}\label{solnoveralgextn} T(E|F)=T(E|F(x)).
\end{equation}

Let $r\in T(K(x)|F).$ Since $F(x)$ is a finite Galois extension of $F,$ so is $K(x)$ over $K.$ Then the ordinary Galois group $Aut(K(x)|K)$ equals the differential Galois group $\gal(K(x)|K)$ (\cite{Put-Singer2003}, Exercise 1.24). Let $\{r=r_1, r_2,\cdots,$  $r_m\}=\{\s(r)\ | \ \s\in \gal(K(x)|K)\}.$  Then for $\s\in \gal(K(x)|K),$   we have $\s(r)\in T(K(x)|F)$ and thus $r_i\in T(K(x)|F)$ for all $i.$ The coefficients of the monic  irreducible polynomial of $r$ over $K$ are symmetric polynomials in $r_1,\cdots, r_m.$ Therefore the coefficients of the irreducible polynomial belong to \[T(K(x)|F)\cap K=T(K|F).\] This shows that $T(K(x)|F)$ is an integral extension of $T(K|F).$ \epf

\section{LIOUVILLIAN PICARD-VESSIOT EXTENSIONS}\label{lpve}

A differential field extension $E$ of $F$ is called a \emph{liouvillian extension} of $F$ if there exists a tower of fields $$E=E_n\subseteq E_{n-1}\supseteq\cdots\supseteq E_0=F$$ such that $E_i=E_{i-1}(t_i)$ and that either $t_i$ is algebraic over $E_{i-1}$ or $t'_i\in E_{i-1}$ or $t_i\neq 0$ and $t'_i/t_i\in E_{i-1}.$ We recall that a Picard-Vessiot extension $E$ of $F$ is called a \emph{liouvillian Picard-Vessiot} extension  if $E$ is a liouvillian extension as well as a Picard-Vessiot extension and that the identity component of the differential Galois group of a  liouvillian Picard-Vessiot extension is solvable.

\begin{theorem}\label{NSliouvillian}
Let $F$ be a differential field with an algebraically closed field of constants and $E$ be a Picard-Vessiot extension of $k.$ Then $E$ is a liouvillian extension of $F$ if and only if $\fieldfrac (T(K|F))=K$ for any differential field $K$ intermediate to $E$ and $F;$ in which case $T(K|F)$ is a finitely generated simple differential $F-$algebra.  
\end{theorem}

\bpf 

\sloppy Let $E$ be a liouvillian extension of $F,$  $K$ be an intermediate differential subfield and $\galsub :=\gal(E|K).$   First we assume that $\gal(E|F)$ is connected. Then we have $T(E|F)\simeq F\otimes_C C[\gal(E|F)]=F[\gal(E|F)].$ Since $\gal(E|F)$ is solvable, from \cite{CPS1977}, Theorem 4.3, we have that  closed subgroups  of $\gal(E|F)$ are  observable. Now from \cite{Hochschild1963}, Theorem $3$, we obtain $\fieldfrac(F[\gal(E|F)]^\galsub)=\fieldfrac(F[\gal(E|F)])^\galsub.$ Thus $\fieldfrac(T(E|F)^\galsub)=\fieldfrac(T(E|F))^\galsub$ and this implies $\fieldfrac(T(K|F))=K.$ Now suppose that $\gal(E|F)$ is not connected. Since $E^{\gal(E|F)^0}=$ $F(x)\subseteq K(x)\subseteq E,$ we have  $\fieldfrac (T(K(x)|F(x))=K(x).$ From Proposition \ref{intext}, $T(K(x)|F(x))=$ $T(K(x)|F)$ and thus $K(x)=\fieldfrac (T(K(x)|F).$  We also know that $T(K(x)|F)$ is an integral extension of $T(K|F).$  Let $S=T(K|F)\setminus \{0\}$. Then  $S^{-1}T(K(x)|F)$ is also an integral extension of $S^{-1}T(K|F).$ Since the latter is a field, so is the former.  However, $K(x)=\fieldfrac (T(K(x)|F))$ is the smallest field containing $T(K(x)|F)$ and thus $ S^{-1}T(K(x)|F)=K(x).$ Now for any $r\in K,$ we have $r=f/g,$ where $f\in T(K(x)|F)$ and $g\in S=T(K|F)\setminus \{0\}.$ Therefore, $f=gr\in T(K(x)|F)\cap K=T(K|F)$ and this proves that $\fieldfrac (T(K|F))=K.$ 

To prove the converse,  we suppose that $E$ is not  a liouvillian extension of $F.$ Let $F(x)=E^{\gal(E|F)^0}$ be the algebraic closure of $F$ in $E.$ Then the identity component $\gal(E|F)^0$  is not solvable and therefore it contains a nontrivial Borel subgroup $\borel.$  Let $K=E^\borel$ and $r\in T(E|F(x))^\borel$.  Since $E$ is a Picard-Vessiot extension of $F(x)$ with Galois group $\gal(E|F)^0,$  the orbit  set $\mathcal{O}_r$ of $r$ under the action of $\gal(E|F)^0 $ is contained in a finite dimensional $C-$vector space which is also $\gal(E|F)^0-$ stable. Moreover, the quotient $\gal(E|F)^0/\borel$ has the structure of a projective variety. Therefore, the induced map $\phi: \gal(E|F)^0/\borel\to \mathcal{O}_r,$ given by $\phi(\bar{\sigma})= \sigma(r)$ for $\sigma\in \gal(E|F)^0,$ is a morphism from a projective variety into some affine space containing  $\mathcal O_r$.  Thus $\phi$  must be a constant. That is, $r\in T(E|F(x))^{\gal(E|F)^0}=F(x)$ and thus $F(x)=T(E|F(x))^{\borel}.$ Note that $$\fieldfrac (T(K|F))=\fieldfrac(T(E|F)^\borel)=\fieldfrac(T(E|F(x))^{\borel})=F(x)\neq K.$$  This proves the converse. 

Next, we shall show that $T(K|F)$ is a finitely generated differential $F-$algebra. First assume that $\gal(E|F)$ is a connected solvable group. Let $\galsub$ be a closed subgroup of $\gal(E|F)$ and $K:=E^{\galsub}.$  We have  $T(E|F)\simeq F\otimes_C C[\gal(E|F)]$ and therefore $$T(K|F)=T(E|F)^\galsub\simeq(F\otimes_C C[\gal(E|F)])^\galsub=F\otimes_C C[\gal(E|F)]^\galsub.$$  Since $\gal(E|F)$ is solvable, the homogeneous space $\gal(E|F)/\galsub$ is affine (\cite{CPS1977}, Theorem 4.3 and Corollary 4.6) and we obtain $C[\gal(E|F)]^\galsub=$  $C[\gal(E|F)/\galsub]$ is a finitely generated $C-$ algebra. This in turn  implies $T(K|F)\simeq F\otimes_C C[\gal(E|F)]^\galsub$ is a finitely generated $F-$algebra. Now assume that only $\gal(E|F)^0$ is solvable.  Let $F(x)=E^{\gal(E|F)^0}$ and observe that $\gal(E|F(x))=\gal(E|F)^0$ is connected.  Then we know $T(K(x)|F(x))$ is a finitely generated $F(x)-$algebra and it follows that $T(K(x)|F(x))$ is a finitely generated $F-$algebra as well.  Since $T(K(x)|F(x))=T(K(x)|F)$ is an integral extension of $T(K|F),$ by Artin-Tate Theorem (\cite{Eisenbud1995}, p.$143$) we obtain that $T(K|F)$ is a finitely generated $F-$algebra.

Now it only remains to show that $T(K|F)$ is a simple differential $F-$algebra. As done earlier, we shall first prove simplicity when $\gal(E|F)$ is connected. Let $I$ be a nonzero  differential ideal of $T(K|F)$ and choose $0\neq y\in I$ so that  $\Lde(y)=0$ for some $\Lde\in F[\partial]$ of smallest positive order $n$.  Since the Galois group is connected, $\Lde=\Lde_{n-1}\Lde_1$ for $\Lde_{n-1}, \Lde_1\in F[\partial]$ of order $n-1$ and $1$ (\cite{Kolchin1948}, p.38). Let $\Lde_1=\partial-a$ for $a\in F$ and observe that $\Lde_1(y)=y'-ay\in I.$ Now since $0=\Lde(y)=\Lde_{n-1}(\Lde_1(y)),$ from the choice of $n$, we obtain that $y'-ay=b\in F.$ Thus $b=\Lde_1(y)\in I.$  If $b\neq 0$ then $I=T(K|F)$. On the other hand if $b=0$ then $y'=ay$ and therefore $(1/y)'=-a(1/y).$ Thus $1/y\in T(K|F)$ and we again obtain  $I=T(K|F).$  This completes the proof when $\gal(E|F)$ is connected.  For an arbitrary liouvillian Picard-Vessiot extension $E$ of $F,$ we have $T(K(x)|F(x))=T(K(x)|F)$ to be a finitely generated simple differential $F-$algebra, where $F(x)$ is the algebraic closure of $F$ in $E$.  Suppose that $T(K|F)$ is not simple and let $I\neq T(K|F)$ be a differential ideal that is maximal among all differential ideals not intersecting $\{1\}$. Then $I$ is known to be a prime ideal.  Let  $I^{e}$ be the extension ideal in $T(K(x)|F).$ It is easy to see that $I^e$ is a differential ideal and therefore $I^e=T(K(x)|F)$. Since $T(K(x)|F)$ is integral over $T(K|F)$ and that $I$ is prime, there must exist a prime ideal of $T(K(x)|F)$ that contracts to $I.$ But any such prime ideal must contain $I^e=T(K(x)|F),$ a contradiction. \epf

\sloppy A liouvillian Picard-Vessiot extension  $E$  of $F$  is known to have the following structure (\cite{Magid1999}, Proposition $6.7$): Let $E^{\gal(E|F)^0}=F(x)$ and $\uni$ be the unipotent radical of $\gal(E|F).$ Then $\uni\subseteq \gal(E|F)^0$ and that  \begin{enumerate}[(a)]\label {unipotent}\item $E=E^\uni(\eta_1,\cdots,\eta_n),$ where $\eta_1,\cdots, \eta_n\in E$ are algebraically independent over $E^\uni$ and $\eta'_i\in E^\uni(\eta_1,\cdots,\eta_{i-1}),$ \item $E^\uni=F(x)(\xi_1,\cdots, \xi_m),$ where $\xi_1,\cdots,\xi_m\in E^\uni$ are algebraically independent over $F(x)$ such that $\xi'_i/\xi_i\in F(x)(\xi_1,\cdots, \xi_{i-1}).$ Moreover, $E^\uni$ is a Picard-Vessiot extension of $F(x)$ with a differential Galois group isomorphic to a maximal torus of $\gal(E|F)^0.$ \end{enumerate}  From the inverse  problem for tori (\cite{Magid1999}, p.99 or \cite{Put-Singer2003}, Exercise 1.41), one can assume further that the $F(x)-$algebraically independent $\xi_i$ are chosen so that $\xi'_i/\xi_i\in F(x)$ (as opposed to $\xi'_i/\xi_i\in F(x)(\xi_1,\cdots, \xi_{i-1})$). Using this description of $E$, in the next corollary, we shall resolve $K$ into a tower of differential fields such that each differential field in the tower  is obtained from its predecessor by  adjoining a solution of a  first order equation of a special kind.

\begin{corollary}\label{structureCLPV}
Let $F$ be a differential field of characteristic zero with an algebraically closed field of constants. Let $E$ be a liouvillian Picard-Vessiot extension of $F$ and  $\gal(E|F)$ be connected. Let $K$ be a differential field intermediate to $E$ and $F.$ Then $K=F(t_1,\cdots, t_n),$ where for each $i,$ $t_i\in T(K|F)$, $t'_i=a_it_i+b_i$ for $a_i\in k$ and $b_i\in F(t_1,\cdots, t_{i-1}).$ Furthermore, if $\gal(E|F)$ is a unipotent algebraic group then each $a_i$ can be taken to be zero and if $\gal(E|F)$ is a torus then each $b_i$ can be taken to be zero.
\end{corollary}

\bpf
To avoid triviality, we shall assume $F\neq K$ and $K\neq E.$ Let $M$ be any differential field such that $F\subseteq M\subsetneqq K.$ We claim that there is a $y\in T(K|F)\setminus M$ such that $y'=ay+b$ for some $a\in F$ and $b\in M$ and that $a$ can be taken to be zero if $\gal$ is unipotent and that $b$ can be taken to be zero if $\gal(E|F)$ is a torus. To prove this claim, we first observe from Theorem \ref{NSliouvillian} that $T(K|F)\setminus M\neq \emptyset$. Choose $y\in T(K|F)\setminus M$ and $\Lde\in F[\partial]$ of smallest positive degree $m$ such that $\Lde(y)=0.$ Since $\gal(E|F)$ is connected,  $\Lde=\Lde_{m-1}\Lde_1,$ where $\Lde_{m-1}, \Lde_1\in F[\partial]$ are of order $m-1$ and $1$ respectively. Let $\Lde_1=\partial-a,$ $a\in F.$ Observe that  $\Lde_1(y)\in T(K|F)$ and $\Lde_{m-1}(\Lde_1(y))=0.$ Therefore, from our choice of $m$, $\Lde_1(y)\in T(M|F)\subset M.$ Thus we have found an element $y\in T(K|F)\setminus M$ such that $y'=ay+b$ where $a\in F$ and $b\in M.$ If $\gal(E|F)$ is unipotent then $E=F(\eta_1,\cdots,\eta_s)$ where $\eta'_i\in k(\eta_1,\cdots,\eta_{i-1})$ and in this case $\Lde$ admits a solution  $\alpha\in F$ (\cite{Ravis2020}, Proposition 2.2). Thus, in this case, we may choose $\Lde_1=\partial-(\alpha'/\alpha)$ and obtain an element $y/\alpha\in T(K|F)$ such that $(y/\alpha)'=b/\alpha\in M.$ Finally suppose that $\gal(E|F)$ is a torus. Then $E=F(\xi_1,\cdots,\xi_s),$
where $\xi'_i/\xi_i\in F$ for each $i.$ If $b\neq 0$ then  apply \cite{Ravis2020}, Proposition 2.2 to the extension $M(\xi_1,\cdots, \xi_s)$ of $M$ with $\Lde(y)=y'-ay=b$  and obtain  $\a\in M$ such that $\a'-a\alpha=b.$ Now $y-\alpha\in T(K|F)\setminus M$ and $(y-\alpha)'/(y-\alpha)=a\in F.$ This proves the claim. Now taking $M=F$ one finds $t_1$ and taking $M=F(t_1,\cdots,t_{i-1}),$ one finds $t_i\in T(K|F)\setminus M$, with the desired properties. Since $K,$ as a field, is finitely generated over $F,$ there must be an $n$ such that $K=F(t_1,\cdots, t_n).$ \epf


\begin{remark} In the above corollary, the  hypothesis that $\gal(E|F)$ is connected allowed us to factor the differential operator $\Lde$ over $F[\partial],$ which was a crucial step in the proof. In fact, the assumption that $\gal(E|F)$ is connected cannot be dropped. For example, consider the liouvillian extension $E=\C(x)(\sqrt{x},e^{\sqrt{x}}),$ where the derivation is $':=d/dx.$ Then $E$ is a liouvillian  Picard-Vessiot extension of $\C(x)$ for the differential equation $$\Lde(Y)=Y''+\frac{1}{2x}Y'-\frac{1}{4x}Y=0.$$ The set  $V:=$span$_\C\{ e^{\sqrt x}, e^{-\sqrt x} \}$ is the set of all solutions of $\Lde(Y)=0$ in $E$.  Since $E$ contains the algebraic extension $\C(x)(\sqrt{x}),$  $\gal(E|F)$ is not connected\footnote{The differential Galois group   $\gal(E|\C(x))$ is isomorphic to $\mathbb G_m \ltimes \mathbb Z_2.$}. One can show that the intermediate differential field  $K:=\C(x)(e^{\sqrt x}+e^{-\sqrt x})$ contains no elements satisfying a first order equation over $\C(x)$ other than the elements of $\C(x)$ itself (\cite{Ravis2020}, p.376).
\end{remark}


\section{INTERMEDIATE DIFFERENTIAL SUBFIELDS OF PICARD-VESSIOT EXTENSIONS}\label{ISPE}
 Let $(\C(x), d/dx)$ be the ordinary differential field of complex rational functions with derivation $':=d/dx.$ Let $E$ be a Picard-Vessiot extension of the Airy differential equation $\Lde(Y):=Y''-xY=0.$ As noted in Section \ref{intro}, for the differential field $K=\C(x)(w),$ where $w$ satisfies the Ricatti equation $w'=x-w^2,$ we have $T(K|\C(x))=\C(x).$  Nonetheless, the differential ring $\C(x)[w]$ is in fact a (finitely generated) simple differential $F-$algebra whose field of fractions is $K$. To see this, it is enough to show  that the differential ring $\C(x)[w]$ is simple.   Suppose that $I$ is a differential ideal of $\C(x)[w].$ Then $I$ must be a principal ideal, say $I=(v)$ for some monic irreducible polynomial $v\in \C(x)[w].$ Write $v=\prod^m_{i=1}w-\alpha_i,$ for distinct algebraic elements $\alpha_i$  of  $\C(x).$ We have $$v'=\sum_j (w'-\alpha'_j)\prod_{i\neq j}(w-\alpha_i).$$ 
Since $v$ divides $v',$ $w-\alpha_i$ must divide $w'-\alpha'_i$  and it follows that $\alpha'_i=x-\alpha^2_i.$ This contradicts the  fact that the Ricatti equation $w'=x-w^2$ has no solutions algebraic over $\C(x).$ Thus $\C(x)[w]$ is a (finitely generated) simple differential $F-$ algebra. This example motivates us to ask whether intermediate differential fields $K$ of arbitrary Picard-Vessiot extensions can be obtained as  field of fractions of some finitely generated simple differential $F-$subalgebras of $K?$  In Theorem \ref{IDSS}, we shall answer this question affirmatively. 

\bprop \label{simdifalg}Let $K$ be a finitely generated
differential field extension of $F$. Then $K$ contains a
finitely generated differential $F-$algebra whose field of
fractions is $K.$ \eprop

\bpf Let $y_1,\cdots, y_{t-1}$ be a transcendence base of $K$ over $F$ and  $F(y_1,\cdots,y_{t-1})[y_t]=K$, for $y_t\in K$ algebraic over $F(y_1,\cdots, y_{t-1}).$ For each $y_i,$ we shall construct a finitely generated differential $F-$algebra $R_i$ whose field of fractions is $F\langle y_i\rangle=F(y_i, y'_i,\cdots).$ Then the smallest $F-$algebra $R$ containing $R_1,\cdots, R_t$ will be a finitely generated differential $F-$algebra whose field of fractions is $K.$ 

Let $y\in \{y_1,\cdots, y_t\}.$ Consider the differential field $F\langle y\rangle$.  Let $n_1$ be the smallest integer so that $y,y',\cdots,$ $y^{(n_1-1)}$ are algebraically independent over $F$ and that $y^{(n_1)}$ be algebraic over the subalgebra  $F[y,y',$ $\cdots, y^{(n_1-1)}]$ of $K$.  Let $$P(X):=\sum^m_{i=0}
a_i X^i\in F[y,y',\cdots, y^{(n_1-1)}][X]$$ be a minimal polynomial of $y^{(n_1)}$ with $a_m\neq 0$.  Now $P(y^{(n_1)})=0$ implies
\begin{equation*}
\sum^m_{i=0} a'_i (y^{(n_1)})^i+\left(\sum^m_{i=0} ia_i
(y^{(n_1)})^{i-1}\right)y^{(n_1+1)}=0.
\end{equation*}
From the minimality of $n_1$, we have $r:=\sum^m_{i=0} ia_i (y^{(n_1)})^{i-1}\neq 0$ and therefore
\begin{equation}\label{higherderivatives}y^{(n_1+1)}=\dfrac{-\sum^m_{i=0} a'_i
(y^{(n_1)})^i}{r_1}\in F[y,y',\cdots,y^{(n_1)}][r^{-1}],\end{equation} Since
$y^{(n_1+1)}\in R_y:=F[y,y',\cdots,y^{(n_1)}, r^{-1}]$, it is
clear for Equation \ref{higherderivatives} that $R_y$ contains all the derivatives of $y$. Also $(1/r)'=-r'/r^2\in R_y$ and thus $R_y$ is a finitely generated differential $F-$algebra whose field of fractions is $F\langle y\rangle$.   
\epf

\bprop\label{Ssimple}

Let $F$ be a differential field of characteristic zero with an algebraically closed
field of constants. Let $\pv$ be a Picard-Vessiot extension of
$F$ and  $S$ be a differential $F-$ subalgebra of $E$ such that $T(E|F)\subseteq S.$ Then $S$ is a simple differential $F-$algebra.
\eprop

\begin{proof} Let $I$ be a nonzero differential ideal of $S$ and $0\neq a\in I.$ Then $a=f/g$ for $f,g\in T(E|F)\subseteq S$ and we have $ga=f\in I^c=T(E|F)\cap I.$ Since the contraction ideal $I^c$ is a differential ideal and that $T(E|F)$ is a simple differential ring, we have $1\in I^c\subseteq I.$ This implies that $S$ is a simple differential ring. 
\end{proof}

\bt\label{IDSS}

Let $F$ be a differential field of characteristic zero with an algebraically closed
field of constants. Let $\pv$ be a Picard-Vessiot extension of
$F$ and let $F\subseteq K\subseteq E$ be an intermediate
differential field. Then $K$ contains a finitely generated simple
differential $F-$algebra whose field of fractions is $K$.

\et

\bpf

Since  Picard-Vessiot extensions are finitely generated field extensions,  we  apply
Proposition \ref{simdifalg} and obtain a finitely generated
differential $F-$algebra $R$, whose field of fractions is $K$. Let
$R:=F[\frac{x_1}{y_1},\frac{x_2}{y_2},\cdots,\frac{x_n}{y_n}],$
where $x_i,y_i\in T(\pv\big/F).$ We shall find an element $r\in R$ so that $R[1/r]$ is a simple differential $F-$algebra and this would complete the proof.

We  first enlarge $T(E|F)$ to a finitely generated simple differential $F-$algebra $S$ so that $R$ becomes a subalgebra of $S.$ To do so, let $S$ be the subring of $E$ generated by $T(\pv\big/F)$ and the set $\{\frac{1}{y_i}\mid
1\leq i\leq n\}$. Since $T(E|F)\subseteq S,$ from Proposition \ref{Ssimple} it follows that $S$ is a finitely generated simple differential algebra.

Next, we shall find a suitable candidate for $r.$ Let $\overline{E}$ be an algebraic closure of $E$ and $\overline{F}$ be the algebraic closure of $F$ in $\overline{E}$. Note that $\overline{F}$ is an algebraically closed field. Let $\R$ and $\Sm$ be the rings generated by $R$ and $S$ over
$\overline{F}$, respectively.  Clearly $R\subset S$, $\R\subseteq \Sm$ and that $\R$ and $\Sm$ are integral extensions of $R$ and $S$ respectively.   The domains  $\R$ and $\Sm$ are finitely generated $\overline{F}-$algebras and therefore they are coordinate rings of  some irreducible affine varieties $X$ and $Y$.  Let $\phi:Y\to X$ be the morphism induced by the inclusion $\R\subseteq \Sm$.  Then $\phi$ is dominant and therefore $\phi(Y)$ must contain an open set $U$ of  $X$. Choose $f\in \R$ so that $X_f:=\{x\in X\ | \ f(x)\neq 0\}\subseteq U$. Since $f$ must be integral over the domain $R$, there is a monic polynomial $P(X)=X^n+r_{n-1}X^{n-1}\cdots+r\in R[X]$ such that $P(f)=0$ and $r\neq 0$. Then $(f^{n-1}+r_{n-1}f^{n-2}+\cdots+r_1)f=-r$ and we have $X_r\subseteq X_f\subseteq U\subseteq \phi(Y)$.  Thus $\phi$ naturally restricts to a surjective morphism from $Y_r$ to $X_r$. Observe that $\variety (I)$ is a non-empty subset of $X_r$  for any proper ideal $I$ of $\R[1/r]$. Since $X_r\subset \phi(Y)$,  we obtain that ${\phi}^{-1}(\variety (I))$ is a non-empty subset of $Y_r$.  Then  ${\phi}^{-1}(\variety (I))\subseteq \variety (I^e)$, where $I^e$ is the extension of $I$ in $\Sm[1/r]$ and therefore $I^e$ is also a proper ideal of $\Sm[1/r]$.  

Now we shall  prove that $R[1/r]$ is a simple differential ring. Suppose that $\id$ is a nonzero proper differential ideal of $R[1/r]$. Since every differential
ideal is contained in a maximal differential ideal and maximal
differential ideals are prime, we have a nonzero prime differential ideal
$\p$ containing $\id$. But $\R[1/r]$ is an integral extension of
$R[1/r]$ and therefore there is a prime ideal $\q$ in
$\R[1/r]$ such that $\q\cap R[1/r]=\p.$ Let $\p^{e}$ be the
extension ideal of $\p$ in $\R[1/r]$. Clearly
$\p^{e}\subseteq\q$ and therefore $\p^e$ is a proper ideal of $\R[1/r]$.  Let $\ib$  be the extension ideal of $\p^e$  in $\Sm[1/r].$ Then from our earlier observation, $\ib$ is a proper   ideal.  Since $\p$ is a differential ideal, so is $\ib$. Now the contraction $\ib^c:=S[1/r]\cap \ib$ must be a proper differential ideal of  $S[1/r].$ Furthermore, $0\neq \p\subseteq R[1/r]\subseteq S[1/r]$ implies $\p\subseteq \ib^c$  and thus $\ib^c$ must be a proper nonzero differential ideal of $S[1/r].$ This contradicts Proposition \ref{Ssimple}.  \epf

 \section*{ACKNOWLEDGEMENT} The authors thank  Chetan Balwe, Andy Magid and Kapil Paranjape for useful discussions.
 
 \bibliographystyle{amsalpha}
\bibliography{LPV-URVRS}


 \end{document}